\documentclass[12pt,letterpaper,oneside,reqno]{amsart}
\usepackage{amsfonts}
\usepackage{amsmath}
\usepackage{amssymb}
\usepackage{amsthm}
\usepackage{float}
\usepackage{mathrsfs}
\usepackage{colonequals}
\usepackage[font=small,labelfont=bf]{caption}
\usepackage[left=1in,right=1in,bottom=1in,top=1in]{geometry}
\usepackage[pdfpagelabels,hyperindex,colorlinks=true,linkcolor=blue,urlcolor=magenta,citecolor=green]{hyperref}
\usepackage{graphicx}
\usepackage{array}
\usepackage{etoolbox}
\apptocmd{\sloppy}{\hbadness 10000\relax}{}{}
\usepackage{derivative}

\newcommand \coeffA [3][A] {{\mathbf{#1}} \sb{#2,#3}}
\newcommand \polynomialP [4][P] {#1 ( #2, #3, #4 )}

\newcommand \derivative [2] {\frac{d}{d #2} #1}                              
\newcommand \qderivative [1] {D_{q} #1}                                      
\newcommand \nqderivative [1] {D_{n,q} #1}                                   
\newcommand \qpowerDerivative [1] {\mathcal{D}_q #1}                         
\newcommand \finiteDifference [1] {\Delta #1}                                
\newcommand \pTsDerivative [2] {\frac{\partial #1}{\Delta #2}}               

\newtheorem{thm}{Theorem}[section]
\newtheorem{cor}[thm]{Corollary}

\newtheorem{examp}[thm]{Example}

\numberwithin{equation}{section}

\title[A study on partial dynamic equation on time scales]
{A study on partial dynamic equation on time scales involving derivatives of polynomials}
\author[Petro Kolosov]{Petro Kolosov}
\email{kolosovp94@gmail.com}
\keywords{
    Dynamic equations on time scales,
    Partial differential equations on time scales,
    Partial dynamic equations on time scales,
    Partial differentiation on time scales,
    Dynamical systems
}
\urladdr{https://kolosovpetro.github.io}
\subjclass[2010]{26E70, 05A30}
\date{\today}
\hypersetup{
    pdftitle={A study on partial dynamic equation on time scales involving derivatives of polynomials},
    pdfsubject={
        Dynamic equations on time scales,
        Partial differential equations on time scales,
        Partial dynamic equations on time scales,
        Partial differentiation on time scales,
        Time-scale calculus,
        Dynamical systems,
        Quantum calculus,
        Differential equations,
        Calculus on time scales,
        Differential calculus,
        Partial differential equations,
        Binomial theorem,
        Discrete convolution,
        Mathematics,
        Classical Analysis and ODEs,
        Analysis of PDEs
    },
    pdfauthor={Petro Kolosov},
    pdfkeywords={
        Dynamic equations on time scales,
        Partial differential equations on time scales,
        Partial dynamic equations on time scales,
        Partial differentiation on time scales,
        Time-scale calculus,
        Dynamical systems,
        Quantum calculus,
        q-calculus,
        q-derivative,
        Differential equations,
        Calculus on time scales,
        Differential calculus,
        Partial differential equations,
        Binomial theorem,
        Binomial coefficients,
        Polynomials,
        Discrete convolution
    }
}
\begin{document}
    \begin{abstract}
        Let $P(m,b,x)$ be a $2m+1$-degree polynomial in $x,b$.
Let be a two-dimensional timescale
$\Lambda^2 = \mathbb{T}_1 \times \mathbb{T}_2 = \{t=(x, b) \colon \; x\in\mathbb{T}_1, \; b\in\mathbb{T}_2 \}$
such that $\mathbb{T}_1 = \mathbb{T}_2$.
In this manuscript we derive and discuss an identity that connects the timescale derivative of odd-power polynomial
with partial derivatives of polynomial $P(m,b,x)$ evaluated in particular points.
For every $t\in\mathbb{T}_1$ and $(x,b) \in \Lambda^2$
\[
    \frac{\Delta x^{2m+1}}{\Delta x}(t) =
    \frac{\partial P(m,b,x)}{\Delta x} (m, \sigma(t), t) +
    \frac{\partial P(m,b,x)}{\Delta b} (m, t, t)
\]
such that $\sigma(t) > t$ is forward jump operator.
In addition, we discuss various derivative operators in context of partial cases of above equation,
we show finite difference, classical derivative, $q-$derivative, $q-$power derivative on behalf of it.

    \end{abstract}

    \maketitle

    \tableofcontents

    \section{Definitions} \label{sec:definitions}
    We now set the following notation such that remains fixed for the remainder of this manuscript
\begin{itemize}
    \setlength\itemsep{1.6em}
    \item Let be a function $f\colon \mathbb{T} \to \mathbb{R}$ and $t\in\mathbb{T}^{\kappa}$ then $f^{\Delta}(t)$
    is delta timescale derivative~\cite{Bohner2001DynamicEO}
    \begin{align*}
        f^{\Delta} (t) = \frac{f(\sigma(t)) - f(t)}{\sigma(t) - t}
    \end{align*}
    where $\sigma(t) - t \neq 0$ and $\sigma(t) > t$ is forward jump operator.

    \item $\dfrac{\partial f(t_1,\ldots,t_n)}{\Delta_i t_i}$ is the delta partial derivative
    of $f\colon \Lambda^n \to \mathbb{R}$ on $n$-dimensional timescale
    $\Lambda^n$ defined via the limit~\cite{bohner2004partial, ahlbrandt2002partial,JACKSON2006391}
    \begin{align*}
        \dfrac{\partial f(t_1,\ldots,t_n)}{\Delta_i t_i} = \lim \limits_{\substack{s_i \to t_i}}
        \frac{
            f(t_1, \ldots, t_{i-1}, \sigma_i(t_i), t_{t+1}, \ldots, t_n)
            - f(t_1, \ldots, t_{i-1}, s_i, t_{t+1}, \ldots, t_n)
        }{\sigma_i(t_i) - s_i}
    \end{align*}
    where $\sigma_i(t_i) > t_i$ and $\sigma_i(t_i) - s_i \neq 0$.

    \item $\qderivative{f(x)}$ is $q-$derivative~\cite{jackson_1909,ernst2000history,ernst2008different,kac2001quantum}
    \begin{align*}
        \qderivative{f(x)} = \frac{f(qx)-f(x)}{qx-x}
    \end{align*}
    where $x\neq 0, \; x\in\mathbb{R}, \; q\in\mathbb{R}$.

    \item $\nqderivative{f(t)}$ is $q-$power derivative~\cite{aldwoah2011power}
    \begin{align*}
        \nqderivative{f(t)} = \frac{f(qt^n) - f(t)}{qt^n - t}
    \end{align*}
    where $qt^n - t \neq 0$ and $n$ is odd positive integer and $0 < q < 1$.

    \item $\qpowerDerivative{f(x)}$ is $q-$power derivative
    \begin{align*}
        \qpowerDerivative{f(x)} = \frac{f(x^q)-f(x)}{x^q-x}
    \end{align*}
    where $x^q \neq x, \; x\in\mathbb{R}, \; q\in\mathbb{R}$.

    \item $\polynomialP{m}{b}{x}$ is $2m+1$-degree polynomial in $x,b$
    \begin{equation}
        \polynomialP{m}{b}{x} = \sum_{k=0}^{b-1} \sum_{r=0}^{m} \coeffA{m}{r} k^r(x-k)^r
        \label{eq:polynomial_p}
    \end{equation}
    where $\coeffA{m}{r}$ is a real coefficient defined recursively, see~\cite{kolosov2016link}.

    \item $\mathbb{Z}$ is an integer timescale such that $\sigma(t) = t+1$.

    \item $\mathbb{R}$ is a real timescale such that $\sigma(t) = t+\Delta t$ where $\Delta t \to 0$.

    \item $q^\mathbb{R}$ is a quantum timescale such that $\sigma(t) = qt$,
    see~\cite[p. 18]{Bohner2001DynamicEO}.

    \item $\mathbb{R}^q$ is a quantum power timescale such that $\sigma(t) = t^q$.

    \item $q^{\mathbb{R}^n}$ is a pure quantum power timescale
    such that $\sigma(t) = qt^n > t, \; 0<q<1$ where $n$ is a positive
    odd integer~\cite{aldwoah2011power}.
\end{itemize}

    \section{Introduction} \label{sec:introduction}
    Time-scale calculus is quite graceful generalization and unification of the theory of differential equations.
Firstly being introduced by Hilger~\cite{hilger1988ein} in his Ph.D thesis in 1988 and thereafter greatly
extended by Bohner and Peterson~\cite{Bohner2001DynamicEO} in 2001, the calculus on time scales became a sharp tool
in the world on differential equations.
Various derivative operators like classical derivative $\derivative{f(x)}{x}$, $q-$derivative $\qderivative{f(x)}$,
$q-$power derivative $\qpowerDerivative{f(x)}$, finite difference $\finiteDifference{f(x)}$ etc,
may be simply expressed in terms of time-scale derivative over particular time scale $\mathbb{T}$.
For instance,
\begin{align*}
    &f^\prime (x) = f^{\Delta} (x), \quad x\in\mathbb{T} = \mathbb{R} \\
    &\Delta f(x) = f^{\Delta} (x), \quad x\in\mathbb{T} = \mathbb{Z} \\
    &\nqderivative{f(x)} = f^{\Delta} (x), \quad x\in\mathbb{T} = q^{\mathbb{R}^n} \\
    &\qderivative{f(x)} = f^{\Delta} (x), \quad x\in\mathbb{T} = q^\mathbb{R} \\
    &\qpowerDerivative{f(x)} = f^{\Delta} (x), \quad x\in\mathbb{T} = \mathbb{R}^q
\end{align*}
In context of Computer Science, namely object oriented programming paradigm, the time scale calculus may be thought
as unified interface of derivative operator.
Furthermore, the idea of time-scale calculus was slightly extended
in~\cite{bayour2017truly,benkhettou2016conformable,caputo2009time,martins2009calculus}.

    \section{Main results} \label{sec:main_results}
    Timescale derivative of odd-powered polynomial $x^{2m+1}$ may be expressed as follows
\begin{thm}
    \label{main_theorem}
    Let $P(m,b,x)$ be a $2m+1$-degree polynomial in $x,b$.
    Let be a two-dimensional timescale
    $\Lambda^2 = \mathbb{T}_1 \times \mathbb{T}_2 = \{t=(x, b) \colon \; x\in\mathbb{T}_1, \; b\in\mathbb{T}_2 \}$
    such that $\mathbb{T}_1 = \mathbb{T}_2$.
    For every $t\in\mathbb{T}_1$ and $(x,b) \in \Lambda^2$
    \begin{align*}
        \frac{\Delta x^{2m+1}}{\Delta x}(t) =
        \frac{\partial P(m,b,x)}{\Delta x} (m, \sigma(t), t) +
        \frac{\partial P(m,b,x)}{\Delta b} (m, t, t)
    \end{align*}
    where
    \begin{itemize}
        \setlength\itemsep{1em}
        \item  $\sigma(t) > t$ -- is forward jump operator

        \item $\frac{\partial \polynomialP{m}{b}{x}}{\Delta x} (m, \sigma(t), t)$ --
        is the value of the partial derivative on time scales of
        $\polynomialP{m}{b}{x}$ with respect to the variable $x$ evaluated in point $(x, b)=  (t, \sigma(t))$

        \item $\frac{\partial \polynomialP{m}{b}{x}}{\Delta b} (m, t, t)$ --
        is the value of the partial derivative on time scales of
        $\polynomialP{m}{b}{x}$ with respect to the variable $b$, evaluated at $(x,b) = (t, t)$
    \end{itemize}
\end{thm}
In simpler words, the theorem ~\ref{main_theorem} says
\begin{center}
    \begin{quotation}
        For every odd-powered polynomial $x^{2m+1}$, the derivative on time scales $\frac{\Delta x^{2m+1}}{\Delta x}$
        evaluated in point $t\in\mathbb{T}_1$ equals to partial derivative on time scales of the polynomial
        $\polynomialP{m}{b}{x}$
        with respect to $x$
        evaluated in point
        $(x,b) = (t, \sigma(t))$
        plus the value of the partial derivative on time scales of the polynomial
        $\polynomialP{m}{b}{x}$
        with respect to $b$,
        evaluated in point
        $(x,b)=(t,t)$.
    \end{quotation}
\end{center}

In its extended form the theorem ~\eqref{main_theorem} is as follows

\begin{align*}
    \frac{\Delta x^{2m+1}}{\Delta x}(t)
    &= \frac{\partial}{\Delta x} \left( \sum_{k=0}^{b-1} \sum_{r=0}^{m} \coeffA{m}{r} k^r(x-k)^r \right) (m, \sigma(t), t) \\
    &+ \frac{\partial}{\Delta b} \left( \sum_{k=0}^{b-1} \sum_{r=0}^{m} \coeffA{m}{r} k^r(x-k)^r \right) (m, t, t)
\end{align*}

    \section{Discussion and examples} \label{sec:discussion_and_examples}
    To understand the nature of the theorem~\ref{main_theorem},
we discuss a few examples involving widely-known time scales
including integer timescale $\mathbb{Z}$,
real timescale $\mathbb{R}$,
quantum timescale $q^{\mathbb{R}}$ and
quantum-power timescale $\mathbb{R}^q$.

\subsection{Time scale of integers $\mathbb{T} = \mathbb{Z} \times \mathbb{Z}$} \label{subsec:time_scale_z}
\begin{cor}
    \label{finite_difference_case}
    (Divided difference.)
    Let be a two-dimensional timescale
    $\Lambda^2 = \mathbb{Z} \times \mathbb{Z}$.
    For every $t\in\mathbb{Z}$ and $x,b\in \Lambda^2$
    \begin{align*}
        \frac{\finiteDifference{x^{2m+1}}}{\Delta x} (t)
        = \pTsDerivative{\polynomialP{m}{b}{x}}{x} (m, \sigma(t), t)
        + \pTsDerivative{\polynomialP{m}{b}{x}}{b}(m, \sigma(t), t)
    \end{align*}
    where $\sigma(t)$ is the forward jump operator defined as $\sigma(t) = t+1$.
\end{cor}
\begin{examp}
    \label{time_scale_z_example_1}
    Let be $t \in \mathbb{Z}, \; x,b \in \Lambda^2 = \mathbb{Z} \times \mathbb{Z}$, let $m=1$ then
    \begin{align*}
        &\pTsDerivative{\polynomialP{1}{b}{x}}{x}                = -3 b + 3 b^2 \\
        &\pTsDerivative{\polynomialP{1}{b}{x}}{b}                = 1 - 6 b^2 + 6 b x
    \end{align*}
    Evaluating in points yields
    \begin{align*}
        &\pTsDerivative{\polynomialP{1}{b}{x}}{x} (1, \sigma(t), t) = 3 t + 3 t^2 \\
        &\pTsDerivative{\polynomialP{1}{b}{x}}{b}(1, t,t)           = 1
    \end{align*}
    Summing up previously obtained partial timescale derivatives, we get ordinary finite difference of odd-powered polynomial
    $x^{3}$ evaluated in point $ t\in\mathbb{Z}, \; x,b\in\Lambda^2 = \mathbb{Z} \times \mathbb{Z}$
    \begin{align*}
        \finiteDifference{x^{3}}(t)
        = \pTsDerivative{\polynomialP{1}{b}{x}}{x} (1, \sigma(t), t)
        + \pTsDerivative{\polynomialP{1}{b}{x}}{b}(1,t,t)
        = 3 t + 3 t^2 + 1
    \end{align*}
\end{examp}
\begin{examp}
    \label{time_scale_z_example_2}
    Let be $t\in\mathbb{Z}, \;x,b\in\Lambda^2 = \mathbb{Z} \times \mathbb{Z}$, let $m=2$ then
    \begin{align*}
        &\pTsDerivative{\polynomialP{2}{b}{x}}{x}                = 5 b - 30 b^2 + 40 b^3 - 15 b^4 + 10 b x - 30 b^2 x + 20 b^3 x \\
        &\pTsDerivative{\polynomialP{2}{b}{x}}{b}                = 1 + 30 b^4 - 60 b^3 x + 30 b^2 x^2
    \end{align*}
    Evaluating in points yields
    \begin{align*}
        &\pTsDerivative{\polynomialP{2}{b}{x}}{x} (1, \sigma(t), t) = 5 t + 10 t^2 + 10 t^3 + 5 t^4 \\
        &\pTsDerivative{\polynomialP{2}{b}{x}}{b} (1, t, t)         = 1
    \end{align*}
    Summing up previously obtained partial timescale derivatives, we get time ordinary finite difference of odd-powered
    polynomial $x^{5}$ and $t\in\mathbb{Z}, \; (x,b) \in\Lambda^2 = \mathbb{Z} \times \mathbb{Z}$
    \begin{align*}
        \finiteDifference{x^{5}} (t)
        = \pTsDerivative{\polynomialP{2}{b}{x}}{x} (1, t, \sigma(t))
        + \pTsDerivative{\polynomialP{2}{b}{x}}{b} (1, t, t)
        = 1 + 5 t + 10 t^2 + 10 t^3 + 5 t^4
    \end{align*}
\end{examp}
\begin{cor}
    \label{time_scale_z_corollary_1}
    For every $t\in\mathbb{Z}, \; (x,b) \in \Lambda^2 = \mathbb{Z} \times \mathbb{Z}$
    \begin{align*}
        \pTsDerivative{\polynomialP{m}{b}{x}}{x} (m, \sigma(t), t)
        = \sum_{r=1}^{2m} \binom{2m+1}{r} t^{r}
    \end{align*}
\end{cor}
\begin{cor}
    \label{time_scale_z_corollary_2}
    For every $t\in\mathbb{Z}, \; (x,b) \in \Lambda^2 = \mathbb{Z} \times \mathbb{Z}$
    \begin{align*}
        \pTsDerivative{\polynomialP{m}{b}{x}}{b} (m, t, t) = 1
    \end{align*}
\end{cor}

\subsection{Time scale of real numbers $\mathbb{T} = \mathbb{R} \times \mathbb{R}$} \label{subsec:time_scale_r}
\begin{cor}
    \label{derivative_case}
    (Classical derivative.)
    Let be a two-dimensional timescale
    $\Lambda^2 = \mathbb{R} \times \mathbb{R} \colonequals \{t=(x, b) \colon \; x\in \mathbb{R}, \; b\in\mathbb{R} \}$.
    For every $t\in\mathbb{R}$ and $(x,b) \in \Lambda^2$
    \begin{align*}
        \odv{x^{2m+1}}{x} (t)
        = \pdv{\polynomialP{m}{b}{x}}{x} (m, \sigma(t), t)
        + \pdv{\polynomialP{m}{b}{x}}{b} (m, t, t)
    \end{align*}
    where $\sigma(t) = t + \Delta t$ such that $ \Delta t \to 0.$
\end{cor}
\begin{examp}
    \label{time_scale_r_example_1}
    Let be $t\in\mathbb{R}, \; (x,b) \in \Lambda^2 = \mathbb{R} \times \mathbb{R}$, let $m=1$ then
    \begin{align*}
        \pdv{\polynomialP{1}{b}{x}}{x} &= -3 b + 3 b^2 \\
        \pdv{\polynomialP{1}{b}{x}}{b} &= 6 b - 6 b^2 - 3 x + 6 b x
    \end{align*}
    Evaluating in points yields
    \begin{align*}
        &\pdv{\polynomialP{1}{b}{x}}{x} (1, \sigma(t), t) = -3t + 3t^2 \\
        &\pdv{\polynomialP{1}{b}{x}}{b} (1, t, t) = 3t
    \end{align*}
    Summing up previously obtained partial timescale derivatives, we get an ordinary derivative of odd polynomial
    $x^{3}$ evaluated in point $t \in \mathbb{R}$.
    \begin{align*}
        \odv{x^3}{x} (t)
        = \pdv{\polynomialP{1}{b}{x}}{x} (1, \sigma(t), t)
        + \pdv{\polynomialP{1}{b}{x}}{b} (1, t, t)
        = 3t^2.
    \end{align*}
\end{examp}
\begin{examp}
    \label{time_scale_r_example_2}
    Let be $t\in\mathbb{R}, \; (x,b) \in \Lambda^2 = \mathbb{R} \times \mathbb{R}$, let $m=2$ then
    \begin{align*}
        &\pdv{\polynomialP{2}{b}{x}}{x} = -15 b^2 + 30 b^3 - 15 b^4 + 10 b x - 30 b^2 x + 20 b^3 x, \\
        &\pdv{\polynomialP{2}{b}{x}}{b} = 30 b^2 - 60 b^3 + 30 b^4 - 30 b x + 90 b^2 x - 60 b^3 x + 5 x^2 - 30 b x^2 + 30 b^2 x^2
    \end{align*}
    Evaluation in points yields
    \begin{align*}
        &\pdv{\polynomialP{2}{b}{x}}{x} (2, \sigma(t), t) = -5 t^2 + 5 t^4 \\
        &\pdv{\polynomialP{2}{b}{x}}{b} (2, \sigma(t), t)        = 5 t^2
    \end{align*}
    Summing up previously obtained partial timescale derivatives, we get classical derivative of an odd polynomial
    $x^5$ evaluated in point $t\in\mathbb{R}$
    \begin{align*}
        \odv{x^5}{x} (t)
        = \pdv{\polynomialP{2}{b}{x}}{x} (2, \sigma(t), t)
        + \pdv{\polynomialP{2}{b}{x}}{b} (2, \sigma(t), t)
        = 5t^4.
    \end{align*}
\end{examp}

\subsection{Quantum time scale $\mathbb{T} = q^\mathbb{R} \times q^\mathbb{R}$} \label{subsec:time_scale_qn}
\begin{cor}
    \label{q_derivative_case}
    (Q-derivative~\cite{jackson_1909}.)
    Let be a two-dimensional time scale
    $\Lambda^2 = q^{\mathbb{R}} \times q^{\mathbb{R}} \colonequals \{t=(x,b) \colon \; x\in q^{\mathbb{R}}, \; b \in q^{\mathbb{R}} \}$.
    For every $t\in q^{\mathbb{R}}$ and $(x,b)\in \Lambda^2$
    \begin{align*}
        \qderivative{x^{2m+1}}(t)
        = \pTsDerivative{\polynomialP{m}{b}{x}}{x} (m, \sigma(t), t)
        + \pTsDerivative{\polynomialP{m}{b}{x}}{b} (m, t, t)
    \end{align*}
    where $\sigma(t) = qt, \; q > 1$.
\end{cor}
\begin{examp}
    \label{time_scale_qn_example_1}
    Let be $t\in q^{\mathbb{R}}, \; x,b\in \Lambda^2 = q^{\mathbb{R}} \times q^{\mathbb{R}}$, let $m=1$ then
    \begin{align*}
        &\pTsDerivative{\polynomialP{1}{b}{x}}{x} = -3 b + 3 b^2 \\
        &\pTsDerivative{\polynomialP{1}{b}{x}}{b} = 3 b - 2 b^2 + 3 b q - 2 b^2 q - 2 b^2 q^2 - 3 x + 3 b x + 3 b q x
    \end{align*}
    Evaluating in points yields
    \begin{align*}
        &\pTsDerivative{\polynomialP{1}{b}{x}}{x} (m, \sigma(t), t) = -3 q t + 3 q^2 t^2 \\
        &\pTsDerivative{\polynomialP{1}{b}{x}}{b} (m, t, t) = 3 q t + t^2 + q t^2 - 2 q^2 t^2
    \end{align*}
    Summing up previously obtained partial time-scale derivatives, we get the $q$-derivative of odd-powered polynomial
    $x^{3}$ evaluated in point $t\in q^{\mathbb{R}}$
    \begin{align*}
        \qderivative{x^{3}}(t)
        = \pTsDerivative{\polynomialP{1}{b}{x}}{x} (m, \sigma(t), t)
        + \pTsDerivative{\polynomialP{1}{b}{x}}{b} (m, t, t)
        = t^2 + q t^2 + q^2 t^2.
    \end{align*}
\end{examp}
For every $t\in q^{\mathbb{R}}, \; (x,b) \in \Lambda^2 = q^{\mathbb{R}} \times q^{\mathbb{R}}$
the following polynomial identity holds as $q$ tends to zero
\begin{align*}
    \lim \limits_{q \to 0} \pTsDerivative{\polynomialP{1}{b}{x}}{b} (1, t, t) = t^2
\end{align*}
However, it would be generalized as follows
\begin{cor}
    \label{time_scale_qn_corollary_1}
    For every $t\in q^{\mathbb{R}}, \; (x,b) \in \Lambda^2 = q^{\mathbb{R}} \times q^{\mathbb{R}}$
    \begin{align*}
        \lim \limits_{q \to 0} \pTsDerivative{\polynomialP{m}{b}{x}}{b} (m, t, t) = t^{2m}.
    \end{align*}
\end{cor}
\begin{examp}
    \label{time_scale_qn_example_2}
    Let be $t\in q^{\mathbb{R}}, \; (x,b) \in \Lambda^2 = q^{\mathbb{R}} \times q^{\mathbb{R}}$, let $m=2$ then
    \begin{align*}
        \pTsDerivative{\polynomialP{2}{b}{x}}{x}
        &= -15 b^2 + 30 b^3 - 15 b^4 + 5 b x - 15 b^2 x + 10 b^3 x + 5 b q x - 15 b^2 q x + 10 b^3 q x \\
        \pTsDerivative{\polynomialP{2}{b}{x}}{b}
        &= 10 b^2 - 15 b^3 + 6 b^4 + 10 b^2 q - 15 b^3 q + 6 b^4 q
        + 10 b^2 q^2 - 15 b^3 q^2 + 6 b^4 q^2 - 15 b^3 q^3 \\
        &+ 6 b^4 q^3 + 6 b^4 q^4 - 15 b x + 30 b^2 x - 15 b^3 x - 15 b q x + 30 b^2 q x
        - 15 b^3 q x + 30 b^2 q^2 x \\
        &- 15 b^3 q^2 x - 15 b^3 q^3 x + 5 x^2 - 15 b x^2 + 10 b^2 x^2 - 15 b q x^2 + 10 b^2 q x^2 + 10 b^2 q^2 x^2
    \end{align*}
    Evaluating in points yields
    \begin{align*}
        &\pTsDerivative{\polynomialP{2}{b}{x}}{x} (2, \sigma(t), t)
        = 5 q t^2 - 10 q^2 t^2 - 15 q^2 t^3 + 15 q^3 t^3 + 10 q^3 t^4 - 5 q^4 t^4 \\
        &\pTsDerivative{\polynomialP{2}{b}{x}}{b} (2, t, t)
        = -5 q t^2 + 10 q^2 t^2 + 15 q^2 t^3 - 15 q^3 t^3 + t^4 + q t^4 + q^2 t^4 - 9 q^3 t^4 + 6 q^4 t^4
    \end{align*}
    Summing up previously obtained partial time-scale derivatives, we get the $q-$derivative of odd polynomial
    $x^{5}$ evaluated in point $t\in q^{\mathbb{R}}$
    \begin{align*}
        \qderivative{t^{5}}
        = \pTsDerivative{\polynomialP{2}{b}{x}}{x} (2, \sigma(t), t)
        + \pTsDerivative{\polynomialP{2}{b}{x}}{b} (2, t, t)
        = t^4 + q t^4 + q^2 t^4 + q^3 t^4 + q^4 t^4.
    \end{align*}
\end{examp}

\subsection{Quantum power time scale $\mathbb{T} = \mathbb{R}^q \times \mathbb{R}^q$} \label{subsec:time_scale_nq}
\begin{cor}
    \label{q_power_derivative_case}
    (Q-power derivative~\cite{aldwoah2011power}.)
    Let be a two-dimensional time scale
    $\Lambda^2 = {\mathbb{R}}^{q} \times {\mathbb{R}}^{q} \colonequals \{t=(x,b) \colon \; b\in {\mathbb{R}}^{q}, \; x\in{\mathbb{R}}^{q} \}$.
    For every $t\in {\mathbb{R}}^{q}, \; (x,b) \in\Lambda^2 = {\mathbb{R}}^{q} \times {\mathbb{R}}^{q}$
    \begin{align*}
        \qpowerDerivative{t^{2m+1}}
        = \pTsDerivative{\polynomialP{m}{b}{x}}{x}(m, \sigma(t), t)
        + \pTsDerivative{\polynomialP{m}{b}{x}}{b}(m, t, t)
    \end{align*}
    where the forward jump operator is defined as $\sigma(t) = t^q, \; q > 1$.
\end{cor}
\begin{examp}
    \label{time_scale_nq_example_1}
    Let be $t\in {\mathbb{R}}^{q}, \; (x,b) \in\Lambda^2 = {\mathbb{R}}^{q} \times {\mathbb{R}}^{q}$, let $m=1$ then
    \begin{align*}
        &\pTsDerivative{\polynomialP{1}{b}{x}}{x} = -3 b + 3 b^2 \\
        &\pTsDerivative{\polynomialP{1}{b}{x}}{b} = 3 b - 2 b^2 + 3 b^q - 2 b^{2 q} - 2 b^{1 + q} - 3 x + 3 b x + 3 b^q x
    \end{align*}
    Evaluating in points yields
    \begin{align*}
        &\pTsDerivative{\polynomialP{1}{b}{x}}{x} (1, \sigma(t), t) = -3 t^q + 3 t^{2 q} \\
        &\pTsDerivative{\polynomialP{1}{b}{x}}{b} (1, t, t) = t^2 + 3 t^q - 2 t^{2 q} + t^{1 + q}
    \end{align*}
    Summing up previously obtained partial time-scale derivatives, we get $q-$power derivative of odd polynomial
    $x^{3}$ evaluated in point $t \in {\mathbb{R}}^{q}$
    \begin{align*}
        \qpowerDerivative{t^{3}}
        = \pTsDerivative{\polynomialP{1}{b}{x}}{x} (1, \sigma(t), t)
        + \pTsDerivative{\polynomialP{1}{b}{x}}{b} (1, t, t)
        = t^2 + t^{2 q} + t^{1 + q}.
    \end{align*}
\end{examp}
\begin{examp}
    \label{time_scale_nq_example_2}
    Let be $t\in {\mathbb{R}}^{q}, \; (x,b) \in \Lambda^2 = {\mathbb{R}}^{q} \times {\mathbb{R}}^{q}$, let $m=2$ then
    \begin{align*}
        \pTsDerivative{\polynomialP{2}{b}{x}}{x}
        &= -15 b^2 + 30 b^3 - 15 b^4 + 5 b x - 15 b^2 x + 10 b^3 x + 5 b x^q - 15 b^2 x^q + 10 b^3 x^q \\
        \pTsDerivative{\polynomialP{2}{b}{x}}{b}
        &= 10 b^2 - 15 b^3 + 6 b^4 + 10 b^{2 q} - 15 b^{3 q} + 6 b^{4 q}
        + 10 b^{1 + q} - 15 b^{2 + q} + 6 b^{3 + q} \\
        &- 15 b^{1 + 2 q} + 6 b^{2 + 2 q} + 6 b^{1 + 3 q} - 15 b x + 30 b^2 x - 15 b^3 x
        - 15 b^q x + 30 b^{2 q} x \\
        &- 15 b^{3 q} x + 30 b^{1 + q} x - 15 b^{2 + q} x - 15 b^{1 + 2 q} x + 5 x^2 - 15 b x^2 + 10 b^2 x^2 \\
        &- 15 b^q x^2 + 10 b^{2 q} x^2 + 10 b^{1 + q} x^2
    \end{align*}
    Evaluation in points yields
    \begin{align*}
        &\pTsDerivative{\polynomialP{2}{b}{x}}{x} (2, \sigma(t), t)
        = -10 t^{2 q} + 15 t^{3q} - 5 t^{4q} + 5 t^{1+q} - 15 t^{1+2q} + 10 t^{1 + 3 q} \\
        &\pTsDerivative{\polynomialP{2}{b}{x}}{b} (2, t, t)
        = t^4 + 10 t^{2 q} - 15 t^{3 q} + 6 t^{4 q} - 5 t^{1 + q} + t^{3 + q}
        + 15 t^{1 + 2 q} + t^{2 + 2 q} - 9 t^{1 + 3 q}
    \end{align*}
    Summing up previously obtained partial time-scale derivatives, we get $q-$power derivative of odd-powered polynomial
    $x^5$ evaluated in point $t\in {\mathbb{R}}^{q}$
    \begin{align*}
        \qpowerDerivative{x^{5}}(t)
        = \pTsDerivative{\polynomialP{2}{b}{x}}{x} (m, \sigma(t), t)
        + \pTsDerivative{\polynomialP{2}{b}{x}}{b} (m, t, t)
        = t^4 + t^{4 q} + t^{3 + q} + t^{2 + 2 q} + t^{1 + 3 q}.
    \end{align*}
\end{examp}
Another polynomial identity, that is exponential sum holds
\begin{cor}
    \label{time_scale_nq_corollary_1}
    For every $t\in {\mathbb{R}}^{q}, \; (x,b) \in \Lambda^2 = {\mathbb{R}}^{q} \times {\mathbb{R}}^{q}, \; t\in\mathbb{R}$
    \begin{align*}
        \lim \limits_{q \to 0} \pTsDerivative{\polynomialP{m}{b}{x}}{b} (m, t, t)
        = \sum_{k=0}^{2m} t^k
    \end{align*}
\end{cor}

\subsection{Pure quantum power time scale $\mathbb{T} = q^{\mathbb{R}^n} \times q^{\mathbb{R}^n}$} \label{subsec:pure_quantum_power}
In this subsection we discuss a pure quantum power time scale $q^{\mathbb{R}^j}$ provided by Aldwoah, Malinowska
and Torres in~\cite{aldwoah2011power}, among with the $q-$power derivative operator $\nqderivative{f(t)}$ defined by
\begin{align*}
    \nqderivative{f(t)} = \frac{f(qt^n) - f(t)}{qt^n - t},
\end{align*}
where $n$ is odd positive integer and $0 < q < 1$.
\begin{cor}
(Quantum power derivative~\cite{aldwoah2011power}.)
    Let be a two-dimensional time scale
    $\Lambda^2 = q^{\mathbb{R}^j} \times q^{\mathbb{R}^j}
    \colonequals \{t=(x,b) \colon \; b \in q^{\mathbb{R}^j}, \; x \in q^{\mathbb{R}^j} \}$.
    For every $t\in q^{\mathbb{R}^j}, \; (x,b) \in \Lambda^2 = q^{\mathbb{R}^j} \times q^{\mathbb{R}^j}$
    \begin{align*}
        \nqderivative{x^{2m+1}} (t)
        = \pTsDerivative{\polynomialP{m}{b}{x}}{x} (m, \sigma(t), t)
        + \pTsDerivative{\polynomialP{m}{b}{x}}{b} (m, t, t)
    \end{align*}
    where $\sigma(t) = qt^n, \; \sigma(t) > t$.
\end{cor}
\begin{examp}
    \label{time_scale_pure_quantum_power_example_1}
    Let be $t\in q^{\mathbb{R}^j}, \; (x,b) \in\Lambda^2 = q^{\mathbb{R}^j} \times q^{\mathbb{R}^j}$, let $m=1$ then
    \begin{align*}
        &\pTsDerivative{\polynomialP{1}{b}{x}}{x} = -3 b + 3 b^2 \\
        &\pTsDerivative{\polynomialP{1}{b}{x}}{b} = 3 b - 2 b^2 + 3 b^j q - 2 b^{1 + j} q - 2 b^{2 j} q^2 - 3 x + 3 b x + 3 b^j q x
    \end{align*}
    Evaluating in points yields
    \begin{align*}
        &\pTsDerivative{\polynomialP{1}{b}{x}}{x}(1, \sigma(t), t) = -3 q t^j + 3 q^2 t^{2 j} \\
        &\pTsDerivative{\polynomialP{1}{b}{x}}{b}(1, t, t) = t^2 + 3 q t^j - 2 q^2 t^{2 j} + q t^{1 + j}
    \end{align*}
    Summing up previously obtained partial timescale derivatives, we get $q-$power derivative of odd-powered polynomial
    $x^{3}$ evaluated in point $t\in q^{\mathbb{R}^j}$
    \begin{align*}
        \nqderivative{x^{3}}(t)
        = \pTsDerivative{\polynomialP{1}{b}{x}}{x} (1, \sigma(t), t)
        + \pTsDerivative{\polynomialP{1}{b}{x}}{b} (1, t, t)
        = t^2 + q^2 t^{2 j} + q t^{1 + j}.
    \end{align*}
\end{examp}

Another polynomial identity, that is exponential sum holds
\begin{cor}
    \label{time_scale_pure_quantum_power_corollary_1}
    For every $t\in q^{\mathbb{R}^j}, \; (x,b) \in\Lambda^2 = q^{\mathbb{R}^j} \times q^{\mathbb{R}^j}, \; t\in\mathbb{R}$
    \begin{align*}
        \lim \limits_{j \to 0} \lim \limits_{q \to 1} \pTsDerivative{\polynomialP{m}{b}{x}}{b} (m, t, t)
        = \sum_{k=0}^{2m} t^k
    \end{align*}
\end{cor}

An identity in even polynomials holds too
\begin{cor}
    \label{time_scale_pure_quantum_power_corollary_2}
    For every $t\in q^{\mathbb{R}^j}, \; (x,b) \in\Lambda^2 = q^{\mathbb{R}^j} \times q^{\mathbb{R}^j}$
    \begin{align*}
        \lim \limits_{j \to 0} \lim \limits_{q \to 0} \pTsDerivative{\polynomialP{m}{b}{x}}{b} (m, t, t)
        = t^{2m}
    \end{align*}
\end{cor}

\begin{examp}
    \label{time_scale_pure_quantum_power_example_2}
    Let be $t\in q^{\mathbb{R}^j}, \; (x,b) \in\Lambda^2 = q^{\mathbb{R}^j} \times q^{\mathbb{R}^j}$, let $m=2$ then
    \begin{align*}
        \pTsDerivative{\polynomialP{2}{b}{x}}{x}
        &= -15 b^2 + 30 b^3 - 15 b^4 + 5 b x - 15 b^2 x + 10 b^3 x + 5 b q x^j - 15 b^2 q x^j + 10 b^3 q x^j \\
        \pTsDerivative{\polynomialP{2}{b}{x}}{b}
        &= 10 b^2 - 15 b^3 + 6 b^4 + 10 b^{1 + j} q - 15 b^{2 + j} q + 6 b^{3 + j} q + 10 b^{2 j} q^2 - 15 b^{1 + 2 j} q^2 \\
        &+ 6 b^{2 + 2 j} q^2 - 15 b^{3 j} q^3 + 6 b^{1 + 3 j} q^3 + 6 b^{4 j} q^4 - 15 b x + 30 b^2 x - 15 b^3 x - 15 b^j q x \\
        &+ 30 b^{1 + j} q x - 15 b^{2 + j} q x + 30 b^{2 j} q^2 x - 15 b^{1 + 2 j} q^2 x - 15 b^{3 j} q^3 x + 5 x^2 - 15 b x^2 \\
        &+ 10 b^2 x^2 - 15 b^j q x^2 + 10 b^{1 + j} q x^2 + 10 b^{2 j} q^2 x^2
    \end{align*}
    Evaluation in points yields
    \begin{align*}
        &\pTsDerivative{\polynomialP{2}{b}{x}}{x} (2, \sigma(t), t)
        = -10 q^2 t^{2 j} + 15 q^3 t^{3 j} - 5 q^4 t^{4 j} + 5 q t^{1 + j} - 15 q^2 t^{1 + 2 j} + 10 q^3 t^{1 + 3 j} \\
        &\pTsDerivative{\polynomialP{2}{b}{x}}{b} (2, t, t)
        = t^4 + 10 q^2 t^{2 j} - 15 q^3 t^{3 j} + 6 q^4 t^{4 j} - 5 q t^{1 + j} + q t^{3 + j} + 15 q^2 t^{1 + 2 j}
        + q^2 t^{2 + 2 j} - 9 q^3 t^{1 + 3 j}
    \end{align*}
    Summing up previously obtained partial timescale derivatives, we $q-$power derivative of odd polynomial
    $x^{5}$ evaluated in point $t\in q^{\mathbb{R}^j}$
    \begin{align*}
        \nqderivative{x^{5}} (t)
        = \pTsDerivative{\polynomialP{1}{b}{x}}{x} (2, \sigma(t), t)
        + \pTsDerivative{\polynomialP{1}{b}{x}}{b} (2, t, t)
        = t^4 + q^4 t^{4 j} + q t^{3 + j} + q^2 t^{2 + 2 j} + q^3 t^{1 + 3 j}
    \end{align*}
\end{examp}

    \bibliographystyle{unsrt}
    \bibliography{AStudyOnDynamicEquations}

\begin{thebibliography}{10}

\bibitem{Bohner2001DynamicEO}
{Bohner, Martin and Peterson, Allan}.
\newblock {\em {Dynamic equations on time scales: An introduction with
  applications}}.
\newblock {Springer Science \& Business Media}, 2001.
\newblock \url {https://web.mst.edu/~bohner/sample.pdf}.

\bibitem{bohner2004partial}
{Bohner, Martin and Guseinov, Gusein Sh.}
\newblock {Partial differentiation on time scales}.
\newblock {\em {Dynamic systems and applications}}, 13(3-4):351--379, 2004.
\newblock \url {http://web.mst.edu/~bohner/papers/pdots.pdf}.

\bibitem{ahlbrandt2002partial}
Calvin~D. Ahlbrandt and Christina Morian.
\newblock Partial differential equations on time scales.
\newblock {\em {Journal of Computational and Applied Mathematics}},
  141(1-2):35--55, 2002.
\newblock \url {https://doi.org/10.1016/S0377-0427(01)00434-4}.

\bibitem{JACKSON2006391}
B.~Jackson.
\newblock Partial dynamic equations on time scales.
\newblock {\em Journal of Computational and Applied Mathematics}, 186(2):391 --
  415, 2006.
\newblock \url {https://doi.org/10.1016/j.cam.2005.02.011}.

\bibitem{jackson_1909}
{Jackson, Frederick H.}
\newblock {XI.—on q-functions and a certain difference operator}.
\newblock {\em {Earth and Environmental Science Transactions of the Royal
  Society of Edinburgh}}, 46(2):253--281, 1909.

\bibitem{ernst2000history}
Thomas Ernst.
\newblock {\em The history of q-calculus and a new method}.
\newblock {Citeseer}, 2000.

\bibitem{ernst2008different}
Thomas Ernst.
\newblock {The different tongues of q-calculus}.
\newblock {\em {Proceedings of the Estonian Academy of Sciences}}, 57(2), 2008.

\bibitem{kac2001quantum}
Victor Kac and Pokman Cheung.
\newblock {\em Quantum calculus}.
\newblock {Springer Science \& Business Media}, 2001.

\bibitem{aldwoah2011power}
{Aldwoah, Khaled A. and Malinowska, Agnieszka B. and Torres, Delfim F.M.}
\newblock The power quantum calculus and variational problems.
\newblock {\em {arXiv preprint arXiv:1107.0344}}, 2011.
\newblock \url {https://arxiv.org/abs/1107.0344}.

\bibitem{kolosov2016link}
{Kolosov, Petro}.
\newblock {On the link between binomial theorem and discrete convolution}.
\newblock {\em {arXiv preprint arXiv:1603.02468}}, 2016.
\newblock \url {https://arxiv.org/abs/1603.02468}.

\bibitem{hilger1988ein}
S.~Hilger.
\newblock {\em {E {\i} n Ma {\ss} kettenkalk {\"u} l m {\i} t Anwendung auf
  Zentrumsmann {\i} gfalt {\i} gke {\i} ten}}.
\newblock PhD thesis, {Ph. D. Thesis, Universt {\"a} t W {\"u} rzburg}, 1988.

\bibitem{bayour2017truly}
{Bayour, Benaoumeur and Hammoudi, Ahmed and Torres, Delfim F.M.}
\newblock {A truly conformable calculus on time scales}.
\newblock {\em {arXiv preprint arXiv:1705.08928}}, 2017.
\newblock \url {https://arxiv.org/abs/1705.08928}.

\bibitem{benkhettou2016conformable}
{Benkhettou, Nadia and Hassani, Salima and Torres, Delfim F.M.}
\newblock A conformable fractional calculus on arbitrary time scales.
\newblock {\em {Journal of King Saud University-Science}}, 28(1):93--98, 2016.

\bibitem{caputo2009time}
{Caputo, M. Cristina}.
\newblock {Time scales: from nabla calculus to delta calculus and vice versa
  via duality}.
\newblock {\em {arXiv preprint arXiv:0910.0085}}, 2009.

\bibitem{martins2009calculus}
{Martins, Natalia and Torres, Delfim F.M.}
\newblock {Calculus of variations on time scales with nabla derivatives}.
\newblock {\em {Nonlinear Analysis: Theory, Methods \& Applications}},
  71(12):e763--e773, 2009.

\bibitem{kolosov2022mathematica}
Petro Kolosov.
\newblock {Supplementary Mathematica Programs}.
\newblock 2020.
\newblock \url
  {https://github.com/kolosovpetro/AStudyOnDynamicEquations/tree/develop/mathematica}.

\end{thebibliography}
    \noindent \textbf{Version:} \texttt{1.0.2-tags-v1-0-1.6+tags/v1.0.1.6d1b96a}

    \clearpage

    \section{Addendum 1: Mathematica scripts} \label{sec:mathematica_scripts}
    To fulfill our study, we attach here a link to the set of \emph{Mathematica} programs, designed to verify the results
of current manuscript.
To reach these programs follow the link~\cite{kolosov2022mathematica}.
To reproduce results, proceed as follows:
\begin{itemize}
    \setlength\itemsep{1em}
    \item Time scale of integers $\mathbb{T} = \mathbb{Z} \times \mathbb{Z}$:
    \begin{itemize}
        \setlength\itemsep{0.5em}
        \item Example~\ref{time_scale_z_example_1}:
        Execute the commands of Mathematica package
        \begin{itemize}
            \item Set \texttt{sigma[x\_] := x + 1} in Mathematica package and execute definition.
            \item Execute \texttt{timeScaleDerivativeX[1, x, b]} which produces $-3 b + 3 b^2$.
            \item Execute \texttt{Expand[timeScaleDerivativeX[1, t, sigma[t]]]} which produces $3 t + 3 t^2$.
            \item Execute \texttt{timeScaleDerivativeB[1, x, b]} which produces $1 - 6 b^2 + 6 b x$.
            \item Execute \texttt{timeScaleDerivativeB[1, t, t]} which produces $1$.
            \item Execute \texttt{mainTheorem[1]} which produces $1 + 3 t + 3 t^2$.
        \end{itemize}
        \item Example~\ref{time_scale_z_example_2}:
        Execute the commands of Mathematica package
        \begin{itemize}
            \item Set \texttt{sigma[x\_] := x + 1} in Mathematica package and execute definition.
            \item \texttt{timeScaleDerivativeX[2, x, b]} which produces
            $5 b - 30 b^2 + 40 b^3 - 15 b^4 + 10 b x - 30 b^2 x + 20 b^3 x$.
            \item \texttt{Expand[timeScaleDerivativeX[2, t, sigma[t]]]} which produces $5 t + 10 t^2 + 10 t^3 + 5 t^4$.
            \item \texttt{timeScaleDerivativeB[2, x, b]} which produces $1 + 30 b^4 - 60 b^3 x + 30 b^2 x^2$.
            \item \texttt{timeScaleDerivativeB[2, t, t]} which produces $1$.
            \item \texttt{mainTheorem[2]} which produces $1 + 5 t + 10 t^2 + 10 t^3 + 5 t^4$.
        \end{itemize}
    \end{itemize}
    \item Time scale of real numbers $\mathbb{T} = \mathbb{R} \times \mathbb{R}$:
    \begin{itemize}
        \item Example~\ref{time_scale_r_example_1}:
        Execute the commands of Mathematica package
        \begin{itemize}
            \item Set \texttt{sigma[x\_] := x + Global`dx} in Mathematica package and execute definition.
            \item Execute \texttt{timeScaleDerivativeX[1, x, b]} which produces $-3 b + 3 b^2$.
            \item Execute \texttt{Limit[Expand[timeScaleDerivativeB[1, x, b]], dx -> 0]}
            which produces $6 b - 6 b^2 - 3 x + 6 b x$.
            \item Execute \texttt{timeScaleDerivativeX[1, t, t]} which produces $-3 t + 3 t^2$.
            \item Execute \texttt{Limit[Expand[timeScaleDerivativeB[1, t, t]], dx -> 0]} which produces $3t$.
            \item Execute \texttt{Limit[mainTheorem[1], dx -> 0]} which produces $3t^2$.
        \end{itemize}
        \item Example~\ref{time_scale_r_example_2}:
        Execute the commands of Mathematica package
        \begin{itemize}
            \item Set \texttt{sigma[x\_] := x + Global`dx} in Mathematica package and execute definition.
            \item Execute \texttt{Limit[Expand[timeScaleDerivativeX[2, x, b]], dx -> 0]}
            which produces $-15 b^2 + 30 b^3 - 15 b^4 + 10 b x - 30 b^2 x + 20 b^3 x$.
            \item Execute \texttt{Limit[Expand[timeScaleDerivativeB[2, x, b]], dx -> 0]}
            which produces $30 b^2 - 60 b^3 + 30 b^4 - 30 b x + 90 b^2 x - 60 b^3 x + 5 x^2 -
            30 b x^2 + 30 b^2 x^2$.
            \item Execute \texttt{Limit[Expand[timeScaleDerivativeX[2, t, sigma[t]]], dx -> 0]} which produces $-5 t^2 + 5 t^4$.
            \item Execute \texttt{Limit[Expand[timeScaleDerivativeB[2, t, t]], dx -> 0]} which produces $5t^2$.
            \item Execute \texttt{Limit[mainTheorem[2], dx -> 0]} which produces $5t^4$.
        \end{itemize}
    \end{itemize}
    \item Quantum time scale $\mathbb{T} = q^\mathbb{R} \times q^\mathbb{R}$:
    \begin{itemize}
        \item Example~\ref{time_scale_qn_example_1}:
        Execute the commands of Mathematica package
        \begin{itemize}
            \item Set \texttt{sigma[x\_] := x * Global`q} in Mathematica package and execute definition.
            \item Execute \texttt{Expand[Simplify[timeScaleDerivativeX[1, x, b]]]}
            which produces $-3 b + 3 b^2$.
            \item Execute \texttt{Expand[Simplify[timeScaleDerivativeB[1, x, b]]]}
            which produces $3 b - 2 b^2 + 3 b q - 2 b^2 q - 2 b^2 q^2 - 3 x + 3 b x + 3 b q x$.
            \item Execute \texttt{Expand[Simplify[timeScaleDerivativeX[1, t, sigma[t]]]]} which produces $-3 q t + 3 q^2 t^2$.
            \item Execute \texttt{Expand[Simplify[timeScaleDerivativeB[1, t, t]]]} which produces $3 q t + t^2 + q t^2 - 2 q^2 t^2$.
            \item Execute \texttt{Expand[Simplify[mainTheorem[1]]]} which produces $t^2 + q t^2 + q^2 t^2$.
        \end{itemize}
        \item Example~\ref{time_scale_qn_example_2}:
        Execute the commands of Mathematica package
        \begin{itemize}
            \item Set \texttt{sigma[x\_] := x * Global`q} in Mathematica package and execute definition.
            \item Execute \texttt{Expand[Simplify[timeScaleDerivativeX[2, x, b]]]}
            which produces $-15 b^2 + 30 b^3 - 15 b^4 + 5 b x - 15 b^2 x + 10 b^3 x + 5 b q x -
            15 b^2 q x + 10 b^3 q x$.
            \item Execute \texttt{Expand[Simplify[timeScaleDerivativeB[2, x, b]]]}
            which produces $10 b^2 - 15 b^3 + 6 b^4 + 10 b^2 q - 15 b^3 q + 6 b^4 q +
            10 b^2 q^2 - 15 b^3 q^2 + 6 b^4 q^2 - 15 b^3 q^3 + 6 b^4 q^3 +
            6 b^4 q^4 - 15 b x + 30 b^2 x - 15 b^3 x - 15 b q x + 30 b^2 q x -
            15 b^3 q x + 30 b^2 q^2 x - 15 b^3 q^2 x - 15 b^3 q^3 x + 5 x^2 -
            15 b x^2 + 10 b^2 x^2 - 15 b q x^2 + 10 b^2 q x^2 + 10 b^2 q^2 x^2$.
            \item Execute \texttt{Expand[Simplify[timeScaleDerivativeX[2, t, sigma[t]]]]}
            which produces $5 q t^2 - 10 q^2 t^2 - 15 q^2 t^3 + 15 q^3 t^3 + 10 q^3 t^4 -
            5 q^4 t^4$.
            \item Execute \texttt{Expand[Simplify[timeScaleDerivativeB[2, t, t]]]}
            which produces $-5 q t^2 + 10 q^2 t^2 + 15 q^2 t^3 - 15 q^3 t^3 + t^4 + q t^4 +
            q^2 t^4 - 9 q^3 t^4 + 6 q^4 t^4$.
            \item Execute \texttt{Expand[Simplify[mainTheorem[2]]]} which produces $t^4 + q t^4 + q^2 t^4 + q^3 t^4 + q^4 t^4$.
        \end{itemize}
        \item Corollary~\ref{time_scale_qn_corollary_1}:
        Execute the commands of Mathematica package
        \begin{itemize}
            \item Set \texttt{sigma[x\_] := x * Global`q} in Mathematica package and execute definition.
            \item Execute \texttt{Limit[Expand[Simplify[timeScaleDerivativeB[m, t, t]]], q -> 0]} for various
            values of \texttt{m}.
        \end{itemize}
    \end{itemize}
    \item Quantum power time scale $\mathbb{T} = \mathbb{R}^q \times \mathbb{R}^q$:
    \begin{itemize}
        \item Example~\ref{time_scale_nq_example_1}:
        Execute the commands of Mathematica package
        \begin{itemize}
            \item Set \texttt{sigma[x\_] := x $\wedge$ Global`q} in Mathematica package and execute definition.
            \item Execute \texttt{Expand[Simplify[timeScaleDerivativeX[1, x, b]]]}
            which produces $-3 b + 3 b^2$.
            \item Execute \texttt{Expand[Simplify[timeScaleDerivativeB[1, x, b]]]}
            which produces $3 b - 2 b^2 + 3 b^q - 2 b^{2 q} - 2 b^{1 + q} - 3 x + 3 b x + 3 b^q x$.
            \item Execute \texttt{Expand[Simplify[timeScaleDerivativeX[1, t, sigma[t]]]]}
            which produces $-3 t^q + 3 t^{2 q}$.
            \item Execute \texttt{Expand[Simplify[timeScaleDerivativeB[1, t, t]]]}
            which produces $t^2 + 3 t^q - 2 t^{2 q} + t^{1 + q}$.
            \item Execute \texttt{Expand[Simplify[mainTheorem[1]]]} which produces $t^2 + t^{2 q} + t^{1 + q}$.
        \end{itemize}
        \item Example~\ref{time_scale_nq_example_2}: Similarly to Example~\ref{time_scale_nq_example_1} with $m=2$.
        \item Corollary~\ref{time_scale_nq_corollary_1}:
        Execute the commands of Mathematica package
        \begin{itemize}
            \item Set \texttt{sigma[x\_] := x $\wedge$ Global`q} in Mathematica package and execute definition.
            \item Execute \texttt{Limit[Expand[Simplify[timeScaleDerivativeB[m, t, t]]], q -> 0]} for various
            values of \texttt{m}.
        \end{itemize}
    \end{itemize}
    \item Pure quantum power time scale $\mathbb{T} = q^{\mathbb{R}^n} \times q^{\mathbb{R}^n}$:
    \begin{itemize}
        \item Example~\ref{time_scale_pure_quantum_power_example_1}:
        Execute the commands of Mathematica package
        \begin{itemize}
            \item Set \texttt{sigma[x\_] := Global`q * x $\wedge$ Global`j} in Mathematica package and execute definition.
            \item Execute \texttt{Expand[Simplify[timeScaleDerivativeX[1, x, b]]]}
            which produces $-3 b + 3 b^2$.
            \item Execute \texttt{Expand[Simplify[timeScaleDerivativeB[1, x, b]]]}
            which produces $3 b - 2 b^2 + 3 b^j q - 2 b^{1 + j} q - 2 b^{2 j} q^2 - 3 x + 3 b x +
            3 b^j q x$.
            \item Execute \texttt{Expand[Simplify[timeScaleDerivativeX[1, t, sigma[t]]]]}
            which produces $-3 q t^j + 3 q^2 t^{2 j}$.
            \item Execute \texttt{Expand[Simplify[timeScaleDerivativeB[1, t, t]]]}
            which produces $t^2 + 3 q t^j - 2 q^2 t^{2 j} + q t^{1 + j}$.
            \item Execute \texttt{Expand[Simplify[mainTheorem[1]]]} which produces $t^2 + q^2 t^{2 j} + q t^{1 + j}$.
        \end{itemize}
        \item Example~\ref{time_scale_pure_quantum_power_example_2}: Similarly
        as Example~\ref{time_scale_pure_quantum_power_example_1} for $m=2$.
        \item Corollary~\ref{time_scale_pure_quantum_power_corollary_1}:
        Execute the commands of Mathematica package
        \begin{itemize}
            \item Set \texttt{sigma[x\_] := Global`q * x $\wedge$ Global`j} in Mathematica package and execute definition.
            \item Execute \texttt{Limit[Limit[Expand[Simplify[timeScaleDerivativeB[m, t, t]]], q -> 1],
                j -> 0]} for various values of \texttt{m}.
        \end{itemize}
        \item Corollary~\ref{time_scale_pure_quantum_power_corollary_2}:
        Execute the commands of Mathematica package
        \begin{itemize}
            \item Set \texttt{sigma[x\_] := Global`q * x $\wedge$ Global`j} in Mathematica package and execute definition.
            \item Execute \texttt{Limit[Limit[Expand[Simplify[timeScaleDerivativeB[5, t, t]]], q -> 0],
                j -> 0]} for various values of \texttt{m}.
        \end{itemize}
    \end{itemize}
\end{itemize}

\end{document}